\documentclass{amsart}

\usepackage{graphicx, epsfig}

\newtheorem{thm}{Theorem}
\newtheorem{lem}{Lemma}
\newtheorem{cl}{Claim}
\newtheorem{defn}{Definition}

\title{An Arctic Circle Theorem for Groves}

\author[T. K. Petersen]{T. Kyle Petersen}
\address{Department of Mathematics, Brandeis University, Waltham, MA, USA, 02454}
\email{tkpeters@brandeis.edu}
\urladdr{http://people.brandeis.edu/\~{}tkpeters}

\author[D. Speyer]{David Speyer}
\address{Department of Mathematics, University of California Berkeley, Berkeley, CA, USA 94720}
\email{speyer@math.berkeley.edu}

\begin{document}

\begin{abstract}
In earlier work, Jockusch, Propp, and Shor proved a theorem
describing the limiting shape of the boundary between the
uniformly tiled corners of a random tiling of an Aztec diamond and
the more unpredictable `temperate zone' in the interior of the
region. The so-called arctic circle theorem made precise a
phenomenon observed in random tilings of large Aztec diamonds.

Here we examine a related combinatorial model called groves.
Created by Carroll and Speyer as combinatorial interpretations for
Laurent polynomials given by the cube recurrence, groves have
observable frozen regions which we describe precisely via
asymptotic analysis of a generating function. Our approach also
provides another way to prove the arctic circle theorem for Aztec
diamonds.
\end{abstract}

\maketitle

\section{Introduction}
Groves came into existence as combinatorial interpretations of the
rational functions generated by the \emph{cube recurrence}:
\[f_{i,j,k}f_{i-1,j-1,k-1} =
f_{i-1,j,k}f_{i,j-1,k-1}+f_{i,j-1,k}f_{i-1,j,k-1} +
f_{i,j,k-1}f_{i-1,j-1,k}, \]  where some initial functions are
specified.  Typically, $f_{i,j,k} := x_{i,j,k}$ for some choice of
$(i,j,k) \in \mathbb{Z}^3$ called the \emph{initial conditions}.
Fomin and Zelevinsky \cite{FZ} were able to show that for initial conditions satisfying some basic requirements, the rational functions generated by the cube recurrence are in fact Laurent polynomials in the
$x_{i,j,k}$. The introduction of groves by Carroll and Speyer
\cite{CS} gave a combinatorial proof of the surprising fact that
each term of these polynomials has coefficient +1. The main
results in this paper only apply to the family of groves on
\emph{standard initial conditions} as described in Section
\ref{sec:init}.\footnote{Herein we will invoke some of the basic
properties of groves without proof. For such arguments, as well as
a general treatment of groves and the cube recurrence, the reader
is referred to \cite{CS}.}

\begin{figure} [h]
\centering
\includegraphics[scale = 1]{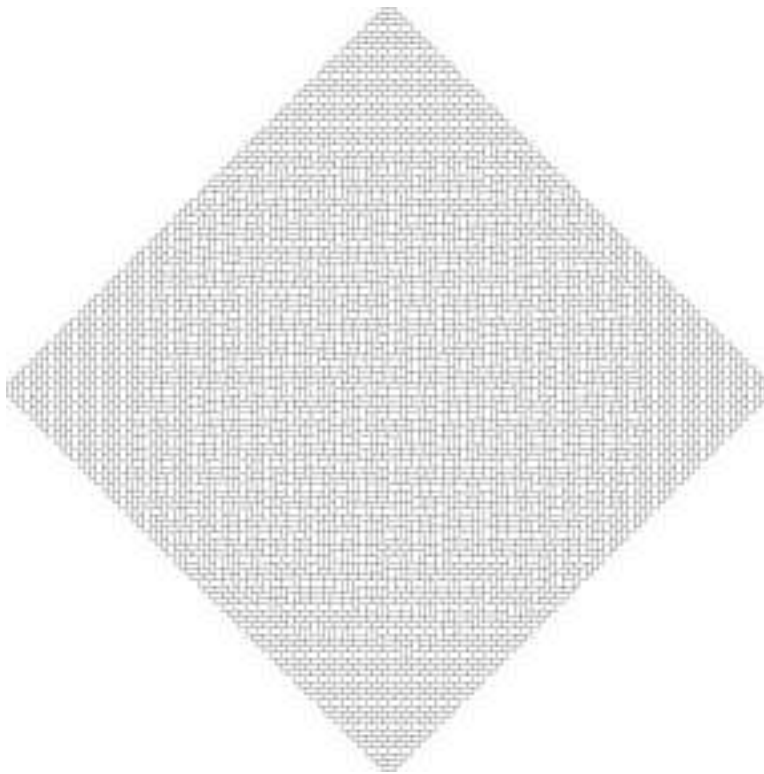}
\caption{A random domino tiling of an Aztec diamond of order 64.}
\label{fig:aztec}
\end{figure}

Before getting into the details of groves, let us first describe
the motivation for this paper: random domino tilings of large
Aztec diamonds.  An \emph{Aztec diamond} of order $n$ consists of
the union of all unit squares with integer vertices contained in
the planar region $\{(x,y) : |x|+|y|\leq n +1\}$.  A \emph{domino
tiling} of an Aztec diamond is an arrangement of $2\times 1$
rectangles, or \emph{dominoes}, that cover the diamond without any
overlapping. A random domino tiling of a large Aztec diamond
consists of two qualitatively different regions.\footnote{By
random we mean selected from the uniform distribution on all
tilings of an Aztec diamond of order $n$, though other probability
distributions may be considered as well.  See \cite{CEP}.}  As
seen in the random tiling in Figure \ref{fig:aztec}, the dominoes
in the corners of the diamond are \emph{frozen} in a brickwork
pattern, whereas the dominoes in the interior have a more random,
\emph{temperate} behavior. It was shown in \cite{JPS} and
\cite{CEP} that asymptotically, the boundary between the frozen
and temperate regions in a random tiling is given by the circle
inscribed in the Aztec diamond. Since everything outside the
circle is expected to be frozen, it is referred to as the
\emph{arctic} circle.

In this paper we shall see that groves on standard initial
conditions exhibit a very similar behavior.  A grove, however, is
not a type of tiling. As the name may suggest, a grove is in fact
a collection of trees. From our point of view, groves are spanning
forests on a finite triangular lattice satisfying certain
connectivity conditions on the boundary. We will show that outside
of the circle inscribed in the triangle, the trees of a large
random grove line up uniformly.

Despite their superficial differences, groves and random domino
tilings of Aztec diamonds are linked by more than their asymptotic
behavior. In fact it seems that their asymptotic behavior is
similar \emph{because} they share a deeper link. The paper of
Carroll and Speyer \cite{CS} establishes that groves are encoded
in the terms of a Laurent polynomial given by the cube recurrence.
There is a more general form of the cube recurrence:
\[f_{i,j,k}f_{i-1,j-1,k-1} = \alpha f_{i-1,j,k}f_{i,j-1,k-1}+
\beta f_{i,j-1,k}f_{i-1,j,k-1} + \gamma f_{i,j,k-1}f_{i-1,j-1,k}\]
where $\alpha, \beta, \gamma$ are constants. If $\alpha= \beta=
\gamma = 1$ we have the original form of the cube recurrence from
whence come groves.  If $\alpha = \beta = 1$ and $\gamma = 0$, we
have (after re-indexing), the \emph{octahedron recurrence}:
\[g_{i,j,n+1}g_{i,j,n-1} =
g_{i-1,j,n}g_{i+1,j,n}+g_{i,j-1,n}g_{i,j+1,n},\] with which we may
encode tilings of Aztec diamonds. In Section \ref{sec:aztec}, we
will show how the polynomial $g_{0,0,n}$ yields all tilings of an
Aztec diamond of order $n$ and we will describe the role that this
recurrence plays in the large scale behavior of such tilings.

While the octahedron recurrence is important to us, it has not
been extant in the study of tilings of Aztec diamonds in the past.
Rather, a local move called \emph{domino shuffling} has been used.
Domino shuffling was introduced in \cite{EKLP} and is generalized
in \cite{propp}. It provides a method for generating tilings of
successively larger Aztec diamonds uniformly at random, and has
been at least implicit in all probabilistic analysis done to date.
Section \ref{sec:shuf} will introduce an analogous local move for
groves that we call \emph{grove shuffling}.  Like domino
shuffling, it will be key to our analysis.

For each of the two models discussed we have a global perspective
and a local perspective. Laurent polynomials tell the global
story: all groves are encapsulated in $f_{0,0,0}$ (from the cube
recurrence), all tilings in $g_{0,0,n}$ (from the octahedron
recurrence). A specified shuffling algorithm tells the local
story. In this paper we combine these two points of view to build
generating functions (for tilings of Aztec diamonds as well as for
groves), with which we can study asymptotic behavior.

\subsection{Groves on standard initial conditions}\label{sec:init}

\begin{figure} [h]
\centering
\includegraphics[scale = .4]{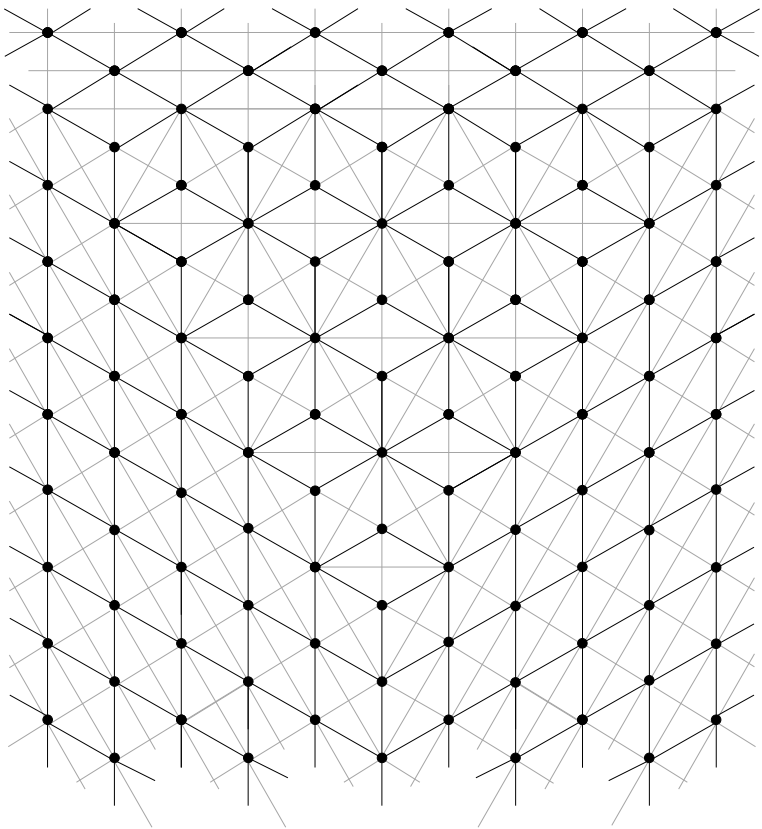}
\caption{A portion of $\mathcal{G}(5)$.}
\label{fig:inits1}
\end{figure}

The standard initial conditions of order $n$ specify a vertex set
$\mathcal{I}(n) = \mathcal{C}(n) \cup \mathcal{B}(n)$ where
$\mathcal{C}(n) = \{ (i,j,k) \in \mathbb{Z}^3 \mid -n-1 \leq i+j+k
\leq -n+1, i,j,k \leq 0\}$ and $\mathcal{B}(n) = \{(i,j,k) \in
\mathbb{Z}^3 \mid i+j+k < -n-1;$ $i,j,k \leq 0;$ and $i,j, \mbox{
or } k = 0\}$.  We draw its projection onto the plane
$\mathbb{R}^{3}/(1,1,1)$ as shown in Figure \ref{fig:inits1} for
the case $n=5$.  One way to generate all groves of order $n$ is to set
$f_{i,j,k} := x_{i,j,k}$ for all $(i,j,k) \in \mathcal{I}(n)$, and
compute $f_{0,0,0}$. Each term in the resulting Laurent polynomial
defines a grove as follows. Let $\mathcal{G}(n)$ be the graph on
the vertex set $\mathcal{I}(n)$ where vertex $(i,j,k)$ has as its
neighbors the vertices $\mathcal{I}(n) \cap \{(i\pm 1,j\pm 1,k),
(i\pm 1,j,k\pm 1), (i,j\pm 1, k\pm 1)\}$. Pictorially, edges of
$\mathcal{G}(n)$ connect vertices that lie diagonally across a
rhombus. In Figure \ref{fig:inits1} the graph $\mathcal{G}(5)$ is made up of the lighter edges and the dark vertices.

\begin{figure} [h]
\centering
\includegraphics[scale = .5]{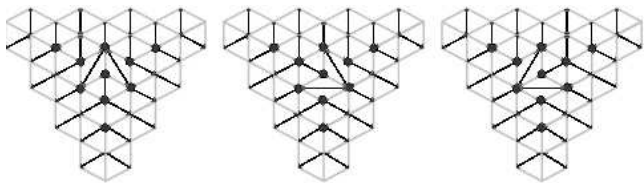}
\caption{The three groves of order 2.} \label{fig:grove1}
\end{figure}

As established in \cite{CS}, the terms in $f_{0,0,0}$ are Laurent
monomials of the form
\[m(g) = \displaystyle \prod_{(i,j,k)\in
\mathcal{I}(n)}x_{i,j,k}^{\deg(i,j,k) - 2},\] where $\deg(i,j,k)
\in \{ 1,2,\ldots,6\}$ is the number of edges connected to vertex
$(i,j,k)$. We have the following
\begin{defn} The \emph{grove} $g$ defined by $m(g)$ is
the unique subgraph of $\mathcal{G}(n)$ containing no crossing
edges such that vertex $(i,j,k)$ in $\mathcal{I}(n)$ has exactly
$\deg(i,j,k)$ incident edges.\end{defn} The uniqueness of the
grove determined by each monomial is a consequence of Theorem 3 in
\cite{CS}. For example, $f_{0,0,0}$ on $\mathcal{I}(2)$ is
\[\frac{x_{-1,-1,0}x_{0,0,-1}}{x_{-1,-1,-1}}+ \frac{x_{-1,0,-1}x_{0,-1,0}}{x_{-1,-1,-1}} + \frac{x_{0,-1,-1}x_{-1,0,0}}{x_{-1,-1,-1}},\]
and the corresponding groves are shown in Figure \ref{fig:grove1}.

\begin{figure} [h]
\centering
\includegraphics[scale=.6]{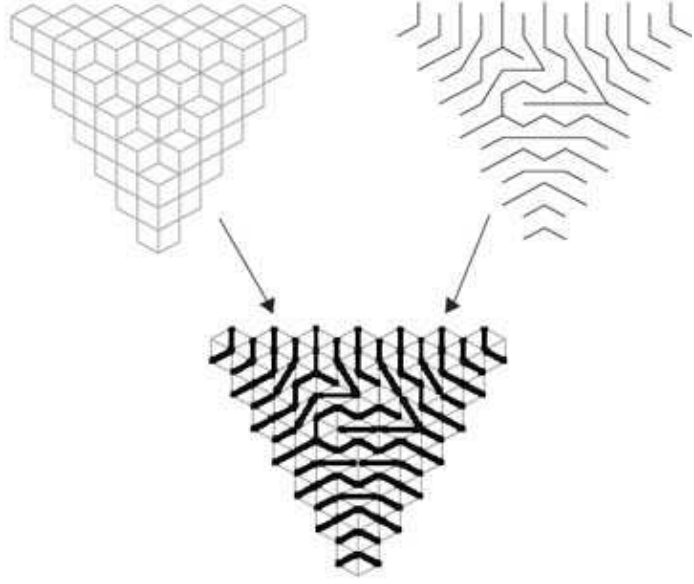}
\caption{A grove $g$ of order 5, superimposed on a picture of
$\mathcal{I}(5)$.}
\label{fig:grove4pic}
\end{figure}

For a more interesting example, one term of $f_{0,0,0}$ on
$\mathcal{I}(5)$ is
\[\frac{x_{-3,0,-2}x_{-2,-1,-1}x_{-1,-3,0}x_{0,-2,-2}}{x_{-3,-1,-2}x_{-2,-3,-1}x_{-1,-2,-2}}.\]
Its corresponding grove, $g$, is shown in Figure
\ref{fig:grove4pic}. We can observe some connectivity properties
of this grove that in fact hold for all groves. Every vertex on the boundary of $\mathcal{C}(n)$
(where cubes have been pushed down) is connected to another vertex
on the boundary of $\mathcal{C}(n)$ if and only if those vertices
are equidistant to the nearest corner (i.e. where two coordinates
are zero) of the grove. Groves are acyclic---every connected component of a grove is a tree. Lastly, each
grove spans $\mathcal{I}(n)$. These connectivity properties are in fact what distinguish groves from arbitrary subgraphs of $\mathcal{G}(n)$, and so give us a combinatorial definition of groves.

Within a grove notice that there are two types of edges: \emph{long} edges and
\emph{short} edges, depending on whether the long or short
diagonal of a rhombus is used. For a vertex $v = (i,j,k)$ in
$\mathcal{C}(n)$, we say that $v$ is:
\begin{description}
\item[up] if $i+j+k = -n+1$,

\item[down] if $i+j+k = -n-1$,

\item[flat] if $i+j+k = -n$,

\item[even] if $i+j+k$ is even,

\item[odd] if $i+j+k$ is odd.
\end{description}
Long edges connect flat vertices to flat
vertices, and short edges connect up vertices to down vertices. Even vertices are only connected to even even vertices and odd vertices are only connected to odd vertices. It
is shown in \cite{CS} that every vertex in $\mathcal{B}(n)$ has
degree 2 and only uses its short edges. As a result, there are
only finitely many long edges, and these determine the grove. This
observation leads to a more convenient way of looking at groves.

\subsection{Simplified groves}
We begin by constructing a modified form of the cube recurrence.
Let $a_{i,j}$, $b_{k,j}$, $c_{i,k}$ be \emph{long edge variables}
where $-n = i+j+k$ is fixed. The variable $a_{i,j}$ is the label
for the edge between vertices $(i,j-1,k+1)$ and $(i-1,j,k+1)$,
$b_{k,j}$ is the label for the edge between $(i-1,j,k+1)$ and
$(i,j,k)$, and $c_{i,k}$ is the label for the edge between
$(i,j,k)$ and $(i,j-1,k+1)$. We write a modified form of the cube
recurrence as follows:
\begin{eqnarray}f_{i,j,k}f_{i-1,j-1,k-1} & = &
b_{i,k}c_{i,j}f_{i-1,j,k}f_{i,j-1,k-1} +
c_{i,j}a_{j,k}f_{i,j-1,k}f_{i-1,j,k-1}\nonumber \\ & & +
a_{j,k}b_{i,k}f_{i,j,k-1}f_{i-1,j-1,k}\nonumber\end{eqnarray} As
we said, the long edges determine the grove, so rather than
setting $f_{i,j,k} := x_{i,j,k}$ for $(i,j,k) \in \mathcal{I}(n)$,
we set $f_{i,j,k} := 1$ for $(i,j,k) \in \mathcal{I}(n)$.  Then
$f_{0,0,0}$ is simply a polynomial in the edge variables
$a_{i,j},b_{i,j},c_{i,j}$, where the variables appear with exponent $+1$ or 0, depending on whether the corresponding long edge is present or not. Each term describes a unique grove, and
we still produce every grove. This form of the cube recurrence is
called the \emph{edge variables version}.

\begin{figure} [h]
\centering
\includegraphics[scale=.8]{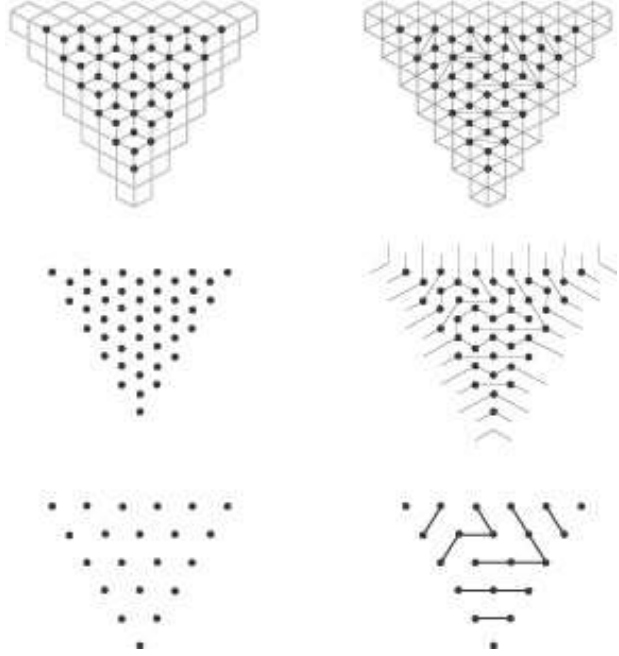}
\caption{On the left: $\mathcal{I}'(5)$ compared to
$\mathcal{I}(5)$. On the right: a standard grove and its
corresponding simplified grove.} \label{fig:simp}
\end{figure}

Taking inspiration from the edge variables version of the cube recurrence, we can draw a simpler
picture of our groves by ignoring all short edges and all of the
vertices incident with them.  In other words, specify a subset of
the standard initial conditions of order $n$, called the
\emph{simplified initial conditions}: $\mathcal{I}'(n) =
\{(i,j,k)\in \mathbb{Z}^3 \mid i + j + k = -n, i,j,k\leq 0\}
\subset \mathcal{I}(n)$. The simplified initial conditions are
just all of the flat vertices. We now represent our groves as
graphs on this vertex set---a triangular lattice shown in Figure
\ref{fig:simp}. Also in Figure \ref{fig:simp} we see the same
grove as in Figure \ref{fig:grove4pic}, but with only the long
edges included. In terms of edge variables, this grove is given by
\[a_{0,0}a_{0,1}a_{0,2}a_{1,0}a_{1,1}a_{2,1}b_{0,0}b_{0,1}c_{0,0}c_{0,1}c_{1,0}c_{2,0}.\]

Another modification of the cube recurrence that we shall like to
use is the \emph{edge-and-face variables version}. In the original
version of the cube recurrence, the variables $x_{i,j,k}$ such
that $i+j+k = -n+1$ were vertex variables.  In the simplified
picture, we call them the \emph{face variables} of order $n$, for
reasons that will become clear. Rather than setting
$f_{i,j,k}:=1$ for all $(i,j,k)$ in $\mathcal{I}(n)$, we give the
face variables their formal weights.  That is, we set
$f_{i,j,k}:=1$ for $(i,j,k) \in \{\,(i,j,k)\in \mathbb{Z}^3 \,|\, -n-1
\leq i+j+k \leq n, i,j,k\leq 0\,\}$ and $f_{i,j,k}:=x_{i,j,k}$ for
$(i,j,k)\in \{\,(i,j,k)\in \mathbb{Z}^3\,|\, i+j+k = -n+1, i,j,k\leq
0\,\}$. Generating $f_{0,0,0}$ using these initial conditions, we
get a Laurent polynomial in the edge and face variables.

The vertices of the simplified initial conditions can be seen as
forming $n(n+1)/2$ downward-pointing equilateral triangles, each
with top-left vertex $(i,j-1,k+1)$, top-right vertex
$(i-1,j,k+1)$, and bottom vertex $(i,j,k)$.  The face variables
then correspond to each of these downward-pointing triangles.  The
triangle with $(i,j,k)$ as its bottom vertex has face variable
$x_{i,j,k+1}$. The exponent of the face variable is $-1, 0,$ or
$1$, corresponding to whether the downward-pointing triangle has,
respectively, two, one, or zero edges present. There can't be
three edges, since that would introduce a cycle and we would no
longer have a forest. Although the face variables don't tell us
anything new about a particular grove, they will be useful later
in deriving probabilities of edges being present in random groves.

\subsection{Grove shuffling}\label{sec:shuf}
We have given one definition for what groves are, and how they may
be generated.  The methods and notation introduced in the previous
section will be very helpful for later proofs.  However, there is
another tool we will like to use; an algorithm called \emph{grove
shuffling} (or \emph{cube-popping} as in \cite{CS}).  Grove
shuffling not only gives a purely combinatorial definition of
groves, but also a method for generating groves of order $n$
uniformly at random. Its inspiration comes from \emph{domino
shuffling}, due to Elkies, Kuperberg, Larsen, and Propp
\cite{EKLP}. The use to which we put grove shuffling is directly
motivated by James Propp and his paper \cite{propp}. For proof
that grove shuffling does indeed give rise to the same objects as
the terms of the Laurent polynomials given by the cube recurrence,
see Carroll and Speyer \cite{CS}. Here we will only include a
description of the algorithm.

\begin{figure} [h]
\centering
\includegraphics[scale=.6]{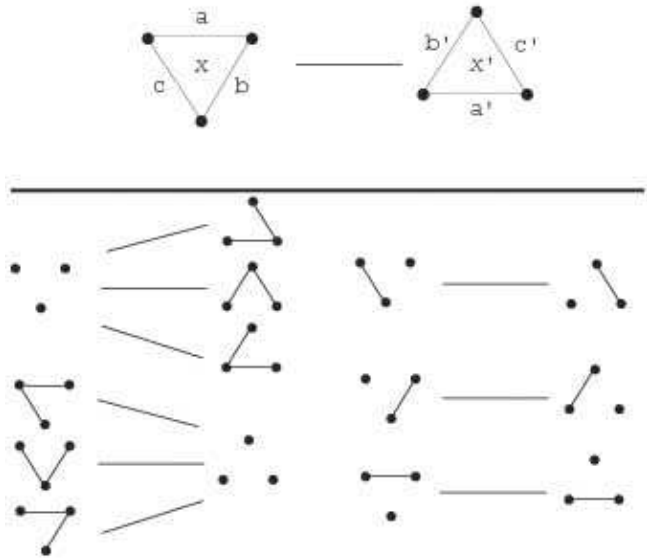}
\caption{Grove shuffling.} \label{fig:shuf}
\end{figure}

\begin{figure} [h]
\centering
\includegraphics[scale=.4]{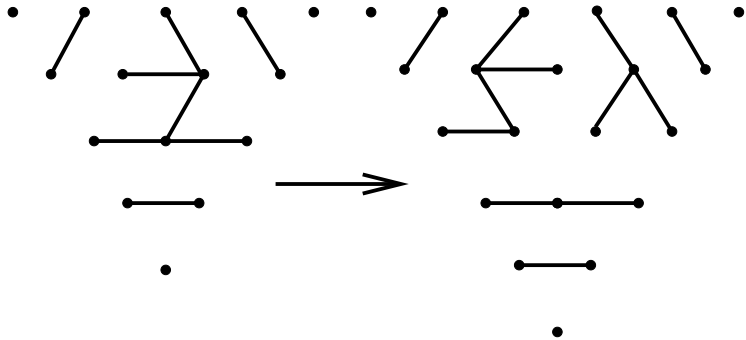}
\caption{A grove of order 4 shuffled into a grove of order 5.} \label{fig:shufexamp}
\end{figure}

Grove shuffling can be thought of as a local move on the
downward-pointing triangles of a simplified grove according to
whether a triangle has zero, one, or two edges present. See Figure
\ref{fig:shuf}. Let $x$ be a generic downward-pointing triangle
with possible edges $a,b,c$ as shown, and let $x'$ be an
upward-pointing triangle, concentric with $x$, with possible edges
$a', b', c'$ as shown.  There are three configurations of $x$ with
two edges: $ab, ac, bc$.  Grove shuffling takes each of these
triangles and replaces them with an upward-pointing triangle $x'$
having none of its possible edges present.  There are three
configurations of $x$ with exactly one edge: $a, b, c$. Each of
these is replaced by the upward-pointing triangle $x'$ with only
the parallel edge: $a', b', c'$, respectively present. Lastly,
there is one configuration of $x$ with none of its possible edges
present. This triangle is replaced with the upward-pointing
triangle $x'$ containing any two of its three possible edges:
$a'b', a'c', b'c'$, chosen randomly with probability 1/3. This
last step is the only random part of the algorithm. After we have
turned every downward-pointing triangle into an upward-pointing
triangle, we add three new vertices to the corners of the grove so
that we may shuffle again. For an example of grove shuffling, see Figure \ref{fig:shufexamp}.\footnote{To see grove shuffling in
action, visit
http://ups.physics.wisc.edu/\~{}hal/SSL/groveshuffler/}

There is a unique grove of order 1. It has one downward-pointing
triangle with zero edges. We now give a purely combinatorial
description of simplified groves on standard initial conditions of
order $n$: they are all the possible results of $n-1$ iterations
of grove shuffling, beginning with the grove of order 1. From
looking at the cube recurrence, it is not hard to show that there
are $3^{\lfloor n^{2}/4\rfloor}$ groves of order $n$. We can now
make the following claim about grove shuffling.

\begin{thm} Beginning with the unique grove of order
one, any grove of order $n$ will be generated after $n-1$
iterations of grove shuffling with probability $1/3^{\lfloor
n^{2}/4\rfloor}$.  In other words, grove shuffling can be used to
generate groves uniformly at random.
\end{thm}

\begin{proof}
Clearly the statement holds for $n$=2.  Suppose that the claim
holds for some $k \geq 1$.  We would like to know the probability
of an arbitrary grove of order $k+1$ being generated.  Fix such a
grove and call it $G(k+1)$. Only a certain subset of the groves of
order $k$ can be shuffled to become $G(k+1)$. Call this set the
shuffling pre-image of $G(k+1)$, denoted $S^{-1}(G(k+1))$. Let
$G(k) \in S^{-1}(G(k+1))$.  Let $a$ be the number of
downward-pointing triangles in $G(k)$ with zero edges, let $b$ be
the number with exactly one edge, and $c$ be the number of
downward-pointing triangles with two edges.

From the rules of grove shuffling, we see that the order of
$S^{-1}(G(k+1))$ is $3^c$. Each pre-image is obtained by making
different choices of the the two edges appearing in each of the
$c$ downward-pointing triangles of $G(k)$.  So since we have
supposed the probability of generating a particular grove of order
$k$ to be uniform, the probability is
\[\frac{3^{c}}{3^{\lfloor\frac{k^2}{4}\rfloor}}\] that after $k$
shuffles we produce a grove in $S^{-1}(G(k+1))$.

Let $S(G(k)) = S(S^{-1}(G(k+1)))$ be the set of groves of order
$k+1$ that can be obtained by shuffling a grove in
$S^{-1}(G(k+1))$. The order of $S(G(k))$ is $3^a$. This is because
in each of the pre-images there are $a$ downward-pointing
triangles with no edges present, and every such triangle can be
shuffled to any of three upward-pointing triangles. Furthermore,
the only edges where the groves of $S^{-1}(G(k+1))$ differ will be
annihilated upon shuffling. So there is a $1/3^a$ chance that one
of the pre-images of $G(k+1)$ will actually shuffle into $G(k+1)$.
Therefore the probability that $k+1$ iterations of grove shuffling
yields $G(k+1)$ is
\[\frac{1}{3^{\lfloor\frac{k^2}{4}\rfloor}}\cdot\frac{1}{3^{a-c}}.\]
Now we claim that $a-c = \lfloor\frac{k+1}{2}\rfloor$.  If so,
then the probability computed above is equal to
\[\frac{1}{3^{\lfloor\frac{(k+1)^2}{4}\rfloor}}\] as desired.

Let us make some basic observations from \cite{CS} or by easy
induction.  First, $a+b+c = k(k+1)/2$; the total number of
downward-pointing triangles in any grove of order $k$.  Secondly,
$b + 2c = \lfloor\frac{k^2}{2}\rfloor$; the total number of edges
in any grove of order $k$.  Then $a-c = k(k+1)/2 -
\lfloor\frac{k^2}{2}\rfloor = \lfloor\frac{k+1}{2}\rfloor$, and
the theorem is proved.
\end{proof}

\subsection{Frozen regions}

We now describe the phenomenon that we analyze in Section
\ref{sec:act}. First we observe that edges are indexed relative to
the corners perpendicular to them, so in fact the edges $a$ and
$a'$ in the description of grove shuffling have the same label: $a
= a' = a_{i,j}$. Horizontal edges are indexed relative to the
bottom corner, and the diagonal edges are indexed relative to the
top-right and top-left corners. In this way we can think of
grove-shuffling as more akin to domino shuffling \cite{propp}.
Rather than replacing edges with parallel edges, we ``slide''
edges toward the corners along perpendicular lines. When a
downward-pointing triangle has two edges, we remove both of those
edges because they ``annihilate'' each other. When a
downward-pointing triangle has no edges, we create two new ones
randomly.

\begin{figure} [h]
\centering
\includegraphics[scale=.7]{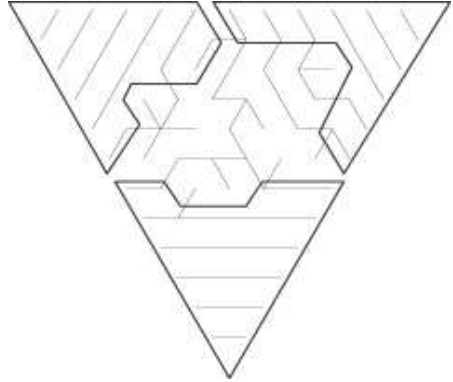}
\caption{Frozen regions of a random grove of order 12.}
\label{fig:frozenexample}
\end{figure}

With this viewpoint, we define an edge to be \emph{frozen} if it
cannot be annihilated under any further iterations of grove
shuffling. Clearly the bottom corner edge, $a_{0,0}$, is frozen
when present. Then the edge $a_{i,j}$ is frozen exactly when the
edges $a_{i',j'}$ are frozen, $i \leq i' \leq 0$, $j \leq j' \leq
0$. Diagonal edges behave similarly. In Figure
\ref{fig:frozenexample} all the highlighted edges are frozen.

\begin{figure} [h]
\centering
\includegraphics[scale=1]{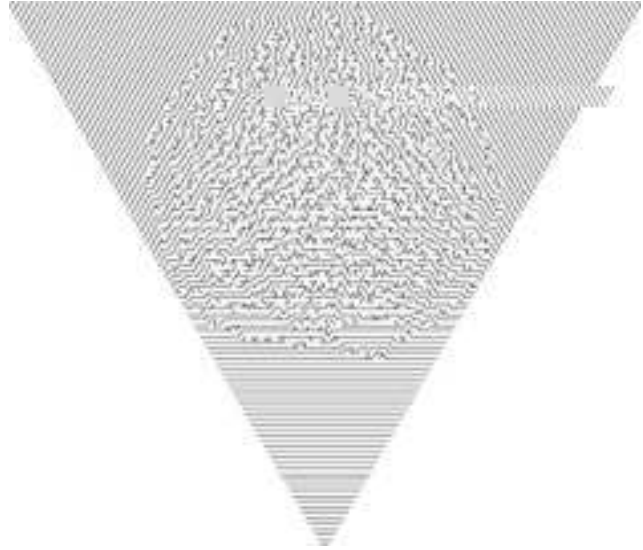}
\caption{A grove on standard initial conditions of order 100.}
\label{fig:ord100}
\end{figure}

We conclude this section by examining a picture of a large random
grove generated by grove shuffling.  In Figure \ref{fig:ord100},
we see that outside of a certain region, all of the edges are
parallel. Moreover, the boundary between the less uniform interior
and the frozen regions in the corners seems to approximate a
circle. Proving that this boundary approaches a circle in the
limit is the main goal of this paper.

\section{The arctic circle theorem}\label{sec:act}

For any $n$, we can scale the initial conditions so that they
resemble an equilateral triangle with sides of length $\sqrt{2}$,
by mapping each vertex $(i,j,k)$ to $(i/n,j/n,k/n)$. The corner
vertices $(-n,0,0), (0,-n,0), (0,0,-n)$ are scaled to $(-1,0,0),
(0,-1,0)$, and $(0,0,-1)$. We will show that outside of the circle
inscribed in this triangle, there is homogeneity of the edges in
an appropriately scaled random grove of order $n$, with
probability approaching 1 as $n \rightarrow \infty$. Specifically,
we will examine the limiting probability of finding a particular
type of edge in a given location outside of the inscribed circle.

\subsection{Edge probabilities}
Let $p_{n}(i,j) = p(i,j,k)$, $k = -n-i-j$, be the probability that
$a_{i,j}(n)$, the horizontal edge on triangle $x_{i,j,k+1}$, is
present in a random grove of order $n$. Similarly define
probabilities $q_{n}(k,i), r_{n}(k,j)$ for the diagonal edges
$b_{k,i}(n)$ and $c_{k,j}(n)$ of the same triangle. Define
$E_{n}(i,j) = E(i,j,k+1)= 1 - p_{n}(i,j) - q_{n}(k,i) -
r_{n}(k,j)$.  The numbers $E_{n}(i,j)$ are analogous to the
\emph{creation rates} discussed in \cite{JPS}, \cite{CEP}, and
\cite{propp}. We will also refer to them as creation rates. As
proven below, we can also realize the number $E_{n}(i,j)$ as the
expected value of the exponent of the face variable $x_{i,j,k+1}$.
We prove the following formula for finding the edge probability
$p_{n}(i,j)$ in terms of creation rates.\footnote{Notice the
similarity between this statement and equation 1.5 of \cite{CEP}.}

\begin{thm}
The horizontal edge probabilities are given recursively by
$p_{n}(i,j) = p_{n-1}(i,j) + \frac{2}{3}E_{n-1}(i,j)$. Thus,
$\displaystyle p_{n}(i,j) =
\frac{2}{3}\sum_{l=1}^{n-1}E_{l}(i,j)$.
\end{thm}

\begin{proof}
We wish to derive a relation between $p_{n}(i,j)$ and
$p_{n-1}(i,j)$. In order to simplify notation, we let:
\begin{displaymath}
\begin{array}{ccc}
p = p_{n-1}(i,j) & q = q_{n-1}(k,i) & r = r_{n-1}(k,j)\\
a = a_{i,j}(n-1) & b = b_{k,i}(n-1) & c = c_{k,j}(n-1)\\
P = p_{n}(i,j) & Q = q_{n}(k,i) & R = r_{n}(k,j)\\
A = a_{i,j}(n) & B = b_{k,i}(n) &C = c_{k,j}(n)
\end{array}
\end{displaymath}
See Figure \ref{fig:prob}.

\begin{figure} [h]
\centering
\includegraphics[scale=1]{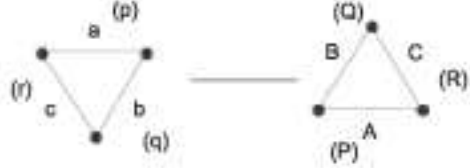}
\caption{Labels for downward- and upward-pointing triangles.}
\label{fig:prob}
\end{figure}

Let $pr(*)$, where $*$ is a subset of $\{a,b,c\}$, be the
probability that a random grove contains that set of edges and not
its compliment. Define $Pr(*)$ similarly.  Some observations that
come directly from grove shuffling:
\begin{itemize}
\item $pr(ab) = pr(ac) = pr(bc)$
\item $pr(abc) = 0$
\item $pr(\emptyset) + pr(a) + pr(b)+pr(c)+pr(ab)+pr(ac)+pr(bc) =
1$
\item $p = pr(a) + pr(ab) + pr(ac)$
\item $q = pr(b) + pr(ab) + pr(bc)$
\item $r = pr(c) + pr(ac) + pr(bc)$
\item $Pr(A) = pr(a)$
\item $Pr(B) = pr(b)$
\item $Pr(C) = pr(c)$
\item $Pr(AB) = Pr(AC) = Pr(BC)= 1/3pr(\emptyset)$
\end{itemize}
We will now deduce $P = p_{n}(i,j)$.
\begin{eqnarray} P
& = & Pr(A) + Pr(AB) + Pr(AC)\nonumber\\
& = & pr(a) + 2/3pr(\emptyset)\nonumber\\
& = & pr(a) + 2/3(1 - pr(a) - pr(b) - pr(c) - pr(ab) -
pr(ac)-pr(bc))\nonumber\\
& = & pr(a) + 2/3(1 - p - q - r + pr(ab)+pr(ac)+pr(bc))\nonumber\\
& = & pr(a) + 2/3(pr(ab) + pr(ac)+pr(bc)) + 2/3(1-p-q-r)\nonumber\\
& = & pr(a) + pr(ab)+pr(ac) + 2/3(1-p-q-r)\nonumber\\
& = & p + 2/3(1-p-q-r) \nonumber
\end{eqnarray}

Let $x= x_{i,j,k+2}$ be the face variable of the downward-pointing
triangle in question.  Notice that
\begin{eqnarray} E(x) & = & \mbox{Expected value of exponent on }x\nonumber\\
& = & 1\cdot pr(\emptyset) + 0\cdot(pr(a) + pr(b) + pr(c)) -
1\cdot(pr(ab)+pr(ac) +pr(bc))\nonumber\\
& = & 1 - pr(a) - pr(b)-pr(c) - 2pr(ab)-2pr(ac)-2pr(bc)\nonumber\\
& = & 1 - p - q-r \nonumber\\
& = & E_{n-1}(i,j).\nonumber
\end{eqnarray}

Therefore, $P = p + 2/3E(x)$.  In the coordinate system, we have
\[p_{n}(i,j) = p_{n-1}(i,j) + 2/3E_{n-1}(i,j) = \frac{2}{3}\sum_{l = 1}^{n-1}E_{l}(i,j)\]
and the theorem is proved.
\end{proof}

\subsection{A generating function}

We now know that to compute the probability of a particular edge
being present in a random grove, it will be enough to compute the
creation rates $E_{l}(i,j)$.  In this section we derive a
generating function for computing these numbers as well as the
related generating function for the horizontal edge probabilities.

Let $\displaystyle F(x,y,z) = \sum_{i,j,k \geq 0}
E(-i,-j,-k)x^{i}y^{j}z^{k}$ be the generating function for the
creation rates. First consider the uniformly weighted version of
the cube recurrence:
\[f_{i,j,k}f_{i-1,j-1,k-1} =
\frac{1}{3}\left(f_{i-1,j,k}f_{i,j-1,k-1}+f_{i,j-1,k}f_{i-1,j,k-1}
+ f_{i,j,k-1}f_{i-1,j-1,k}\right).\] We will return to the
convention of setting $f_{i,j,k} =x_{i,j,k}$ for all $(i,j,k)\in
\mathcal{I}(n)$. Using this recurrence to calculate $f_{0,0,0}$ we
will get each monomial weighted uniformly, so that if we set all
the variables equal to 1, $f_{0,0,0} = 1$. If we want the
expectation of the exponent of the face variable $x =
x_{i_{0},j_{0},k_{0}}$, we need only calculate the derivative of
$f_{0,0,0}$ with respect to this variable, then set all variables
equal to one.  In other words,
\[E(i_{0},j_{0},k_{0}) = \frac{\partial}{\partial x}\Big(f_{0,0,0}\Big)\Big|_{x_{i,j,k} =
1}\] Furthermore, we can calculate the intermediate creation rates
the same way.
\begin{lem}
Fix $x = x_{i_0,j_0,k_0}$ for $i_0+j_0+k_0 = -n+1$. Then for any
$(i',j',k')$ such that $i'+j'+k' = -n'+1$ with $n' < n$, we have
\[\frac{\partial}{\partial
x}\Big(f_{i',j',k'}\Big)\Big|_{x_{i,j,k} = 1} =
E(i_0-i',j_0-j',k_0-k').\]
\end{lem}

\begin{proof}
First, we re-center our initial conditions to make the situation
clear. Introduce the variables $\tilde{x}$ by $x_{i,j,k} \to
\tilde{x}_{i-i',j-j',k-k'}$. In particular, $x_{i_0,j_0,k_0} \to
\tilde{x}_{i_0 -i',j_0-j',k_0 -k'}$. Then we have $f_{i',j',k'}
\to  \tilde{f}_{0,0,0}$, a polynomial generated from a set of
standard initial conditions of order $n-n'-1$. Differentiating
$f_{i',j',k'}$ with respect to $x_{i_0,j_0,k_0}$ and setting all
the variables to one is equivalent to differentiating
$\tilde{f}_{0,0,0}$ with respect to $\tilde{x}_{i_0 -i',j_0-j',k_0
-k'}$ and setting all variables equal to one. The latter action
clearly gives $E(i_0 -i',j_0-j',k_0 -k')$, and the claim is
proved.
\end{proof}

With this in mind, let us differentiate the weighted cube
recurrence with respect to $x=x_{i_0,j_0,k_0}$:
\begin{eqnarray*}
f'_{i,j,k}f_{i-1,j-1,k-1} + f_{i,j,k}f'_{i-1,j-1,k-1} & = &
\frac{1}{3}\Big(f'_{i-1,j,k}f_{i,j-1,k-1} +
f_{i-1,j,k}f'_{i,j-1,k-1} + \\
& & f'_{i,j-1,k}f_{i-1,j,k-1} +
f_{i,j-1,k}f'_{i-1,j,k-1} + \\
& & f'_{i,j,k-1}f_{i-1,j-1,k} + f_{i,j,k-1}f'_{i-1,j-1,k}\Big).
\end{eqnarray*}
Now by setting $x_{i,j,k} = 1$ for all $(i,j,k)$, we get a linear
recurrence for the expectations in question (where $r=i_0-i,
s=j_0-j, t=k_0-k$):
\begin{eqnarray*}
E(r,s,t) + E(r+1,s+1,t+1) & = & \frac{1}{3}\Big(E(r+1,s,t) +
E(r,s+1,t+1) +\\
& & E(r,s+1,t) + E(r+1,s,t+1) +
\\
& & E(r,s,t+1) + E(r+1,s+1,t)\Big).
\end{eqnarray*}
The recurrence holds for any $r,s,t < 0$. Let us also observe some
recurrences near the boundary.
\begin{eqnarray*}
E(r,0,0) &=& \frac{1}{3}E(r+1,0,0)\\
E(0,s,0) &=& \frac{1}{3}E(0,s+1,0)\\
E(0,0,t) &=& \frac{1}{3}E(0,0,t+1)\\
E(r,s,0) &=&
\frac{1}{3}\Big(E(r+1,s,0)+E(r,s+1,0)+E(r+1,s+1,0)\Big)\\
E(r,0,t) &=&
\frac{1}{3}\Big(E(r+1,0,t)+E(r,0,t+1)+E(r+1,0,t+1)\Big)\\
E(0,s,t) &=&
\frac{1}{3}\Big(E(0,s+1,t)+E(0,s,t+1)+E(0,s+1,t+1)\Big)
\end{eqnarray*}
After computing $E(0,0,0) = 1$, we can form the rational
generating function in the variables $x,y,z$:
\begin{eqnarray*}
F(x,y,z) & = & \sum_{i,j,k \geq 0} E(-i,-j,-k)x^{i}y^{j}z^{k}
 \\
& = & \frac{1}{1+xyz-\frac{1}{3}(x+y+z+xy+xz+yz)}
\end{eqnarray*}
Now using the fact that $p(i,j,k) = p(i,j,k+1) + (2/3)E(i,j,k+2)$,
we can derive the formula for the generating function we want:
\begin{eqnarray*}
G(x,y,z) &= & \sum_{i,j,k \geq 0}p(-i,-j,-k)x^{i}y^{j}z^{k} \\
& = & \frac{2z^2F(x,y,z)}{3(1-z)}.
\end{eqnarray*}

\subsection{Asymptotic analysis}

With our generating function in hand, we can prove our main
theorem.  First let us embed a triangle in three-space by $T:=
\{(x,y,z)\in \mathbb{R}^3 \mid x,y,z \leq 0, x+y+z = -1\}$. This
is the triangle that we will scale $\mathcal{I}(n)$ to fit.  A
point $(x,y,z) \in T$ is outside of the inscribed circle (what
will show is the arctic circle) if and only if the angle between
the vector $(x,y,z)$ and vector $(-1,-1,-1)$ is greater than
$\cos^{-1}(\sqrt{2/3})$.

Notice that for any point $(x,y,z)$ outside of the inscribed
circle, we have either $x,y > -1/2$, $x,z > -1/2$, or $y,z >
-1/2$. We call any coordinate with a value strictly greater than
$-1/2$ a \emph{small} coordinate. We now state

\begin{thm}[Weak Arctic Circle]\label{thm:wac}
Let $(x_{0},y_{0},z_{0})$ be a point in $T$ outside of the
inscribed circle for which $z_{0}$ is a small coordinate.  Let
$(i_{n},j_{n},k_{n})$, $i_{n}+j_{n}+k_{n} = -n$, be a sequence of
nonpositive integer triples such that
\[\lim_{n\to\infty}\frac{1}{n}(i_{n},j_{n},k_{n}) = (x_{0},y_{0},z_{0})\]
Then $\displaystyle \lim_{n\to\infty} p(i_{n},j_{n},k_{n}) = 0$.
\end{thm}

In other words, the theorem states that in the upper two regions
of $T$ outside of the arctic circle, the probability of finding a
horizontal edge goes to zero as the order of a (scaled) random
grove goes to infinity.  By symmetry, there can be no diagonal
edges in the lower region, and in order to satisfy the
connectivity properties of groves, all the edges in the lower
region must be horizontal. The following lemma is the heart of the
proof.

\begin{lem}\label{lem1}
Fix a point $(x_{0},y_{0},z_{0})$ in $T$ outside of the inscribed
circle for which $z_{0}$ is a small coordinate. Then there are
real constants $A,B,C$ such that
\[p(-i,-j,-k) = O(e^{-(Ai+Bj+Ck)})\] for all $i,j,k \geq 0$ and $Ax_{0}+By_{0}+Cz_{0} <
0$.
\end{lem}

Proving Lemma \ref{lem1} consumes most of the rest of this
section. Let us suppose the lemma is true and present the proof of
the theorem.
\begin{proof}[Proof of Theorem \ref{thm:wac}]
By Lemma \ref{lem1}, $p(i_{n},j_{n},k_{n}) =
O(e^{Ai_{n}+Bj_{n}+Ck_{n}})$, so we will have that
$p(i_{n},j_{n},k_{n})\to 0$ if $Ai_{n}+Bj_{n}+Ck_{n}\to -\infty$.

Say $Ax_{0}+By_{0}+Cz_{0}=d < 0$ and let \[a_{n} = A\frac{i_n}{n}
+ B\frac{j_n}{n} + C\frac{k_n}{n}.\] Then for any $\epsilon > 0$
there is some $N > 0$ such that for all $n \geq N$, $a_{n} \in
B(d,\epsilon) = \{\,x : |d-x| < \epsilon\,\}$. So if we take $n$
  sufficiently large, $n a_{n} < 0$, and $|d+\epsilon|n < |a_{n}|n$.
Since $|d+\epsilon|n \to \infty$, by comparison we have
\[n a_{n} = Ai_n + Bj_n + Ck_n \to -\infty,\] and the theorem
is proved.
\end{proof}

To facilitate the proof of Lemma \ref{lem1}, we will use the following
claims.

\begin{cl}\label{cl1}
Let $f(x,y,z)$ be an analytic function.  Let $r,s,t$ be positive
real numbers such that $f(x,y,z)\neq 0$ for $|x|\leq r$, $|y|\leq
s$, $|z|\leq t$.  If \[G(x,y,z) = \frac{1}{f(x,y,z)} =
\sum_{i,j,k\geq 0} a_{i,j,k}x^{i}y^{j}z^{k},\] then $a_{i,j,k} =
O(r^{-i}s^{-j}t^{-k})$.
\end{cl}

\begin{proof}[Proof of Claim \ref{cl1}]
Define the loops $\gamma = \{|z|=t\}, \gamma' = \{|y|=s\},$ and
$\gamma'' = \{|x|=r\}$. Then we have
\begin{eqnarray}
a_{i,j,k} & = & \frac{1}{i!j!k!}\frac{\partial}{\partial
x^{i}}\frac{\partial}{\partial y^{j}}\frac{\partial}{\partial
z^{k}}G(x,y,z)\Big|_{(0,0,0)}\nonumber\\
& = & \frac{1}{(2\pi
i)^{3}}\int_{\gamma}\int_{\gamma' }\int_{\gamma''}\frac{G(x,y,z)}{x^{i+1}y^{j+1}z^{k+1}}dxdydz\nonumber\\
& \leq & \frac{M}{(2\pi i)^{3}}\int_{\gamma}\int_{\gamma'
}\int_{\gamma''}\frac{1}{x^{i+1}y^{j+1}z^{k+1}}dxdydz\nonumber\\
& & (\mbox{Since } G(x,y,z) \mbox{ is bounded on the compact
set } \gamma\times \gamma' \times \gamma''.)\nonumber\\
& = & M\frac{1}{r^{i}s^{j}t^{k}}\nonumber
\end{eqnarray}
In other words, $a_{i,j,k} = O(r^{-i}s^{-j}t^{-k})$ and the claim
is proved.
\end{proof}

\begin{cl}\label{cl2}
Let $(x_0,y_0,z_0)$ be as in Lemma \ref{lem1}. Without loss of
generality, say $y_0$ is the second small coordinate. Then there
exists a vector $(A_0,B_0,C_0) \in \mathbb{R}^3$ such that
\begin{itemize}
\item $A_0 x_{0}+B_0 y_{0}+C_0 z_{0} < 0$
\item $A_0 B_0+A_0 C_0 + B_0 C_0 > 0$
\item $B_0, C_0 < 0$.
\end{itemize}
\end{cl}

\begin{proof}[Proof of Claim \ref{cl2}]
Let $\phi$ be the angle between $(x_0,y_0,z_0)$ and $(-1,-1,-1)$.
The assumption that $(x_0,y_0,z_0)$ lies outside the Arctic circle
means $\phi > \sin^{-1} (\sqrt{1/3})$. Take $(A_0,B_0,C_0)$ in the
plane spanned by $(-1,-1,-1)$ and $(x_0,y_0,z_0)$ and lying on the
other side of $(-1,-1,-1)$ from $(x_0,y_0,z_0)$. Let $\theta$ be
the angle between $(A_0,B_0,C_0)$ and $(-1,-1,-1)$, so $\theta +
\phi$ is the angle between $(A_0,B_0,C_0)$ and $(x_0,y_0,z_0)$. If
$\phi = \sin^{-1} (\sqrt{1/3})+\epsilon$, $\epsilon > 0$, we can
choose $\theta = \cos^{-1} (\sqrt{1/3}) -\epsilon/2$ so that
$\theta+\phi =\cos^{-1} \sqrt{1/3}+\sin^{-1} \sqrt{1/3} +
\epsilon/2 > \cos^{-1} \sqrt{1/3}+\sin^{-1} \sqrt{1/3}=\pi/2$.

\begin{figure} [h]
\centering
\includegraphics[scale =.7]{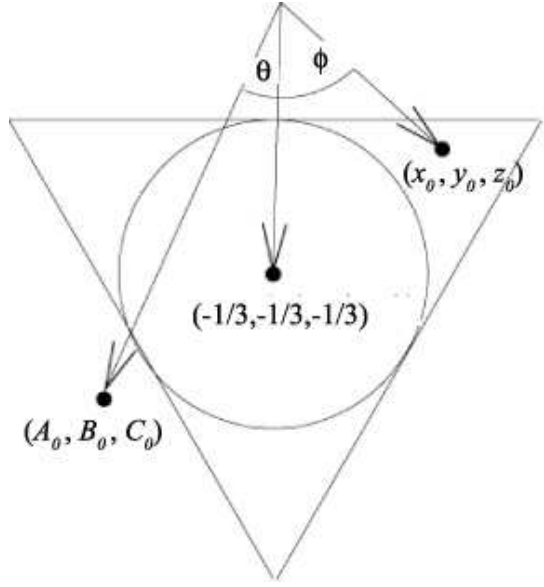}
\caption{The relative positions of $(x_0,y_0,z_0)$ and $(A_0,B_0,C_0)$.}
\label{fig:cone}
\end{figure}

From $\theta< \cos^{-1} \sqrt{1/3}$, we deduce
\begin{eqnarray*}
\frac{(A_0+B_0+C_0)^2}{3(A_0^2+B_0^2+C_0^2)} & = &
\frac{\langle(A_0,B_0,C_0),(-1,-1,-1)\rangle^2}{\|(A_0,B_0,C_0)\|^2\|(-1,-1,-1)\|^2}
\\
 & = & \cos^2 \theta \\
 & > & \frac{1}{3}
\end{eqnarray*}
so $(A_0+B_0+C_0)^2 > A_0^2+B_0^2+C_0^2$ which is equivalent to
$A_0 B_0+A_0 C_0 + B_0 C_0 > 0$. From $\theta+\phi>\pi/2$, we
deduce that $A_0 x_0+B_0 y_0 + C_0 z_0 <0$.

Finally, we must show that we can choose $B_0, C_0 < 0$. Let
$\alpha$ be the distance between $(A_0, B_0, C_0)$ and
$1/3(-1,-1,-1)$. Let $\beta > \sqrt{3}/6$ be the distance between
$(x_0, y_0, z_0)$ and $1/3(-1,-1,-1)$. We have, for
$(A_0,B_0,C_0)$ satisfying the first two conditions and lying in
the plane $x+y+z = -1$, that $\alpha\beta > 1/6$. Now let $\vec{v}
= (x_0,y_0,z_0)+1/3(1,1,1)$ so that $|\vec{v}| = \beta$. Then
choosing
\[(A_0,B_0,C_0) = -\alpha\frac{\vec{v}}{|\vec{v}|} -\frac{1}{3}(1,1,1)
= -\frac{\alpha}{\beta}\vec{v} -\frac{1}{3}(1,1,1)\] has all the
desired properties. The first two are obvious by construction, and
using the facts that $-\alpha < -1/(6\beta)$, $\beta
> \sqrt{3}/6$ and, because $z_0$ is a small coordinate, $-z_0 < 1/2$, we can verify that
\begin{eqnarray*}
C_0 & = & -\frac{\alpha}{\beta}(z_0 + 1/3) - 1/3\\
 & < & \frac{-z_0 - 1/3}{6\beta^{2}} - 1/3 \\
 & < & \frac{1/2 - 1/3}{6\beta^{2}} - 1/3 \\
 & = & \frac{1}{36\beta^2} - 1/3 \\
 & < & 1/3 - 1/3 = 0
\end{eqnarray*}
If $y_0$ is the other small coordinate then an identical
computation shows $B_0 < 0$.
\end{proof}

\begin{cl}\label{cl3}
Let $(A,B,C) \in \mathbb{R}^3$ with $C<0$. Suppose that \[g(x,y,z)
= 1+xyz - (1/3)(x+y+z+xy+xz+yz)\] is not zero for $(x,y,z) \in
[0,e^A] \times [0,e^B] \times [0,e^C]$. Then $g(x,y,z)$ is not
zero for any $(x,y,z) \in \mathbb{C}^3$ with $(x,y,z) \in
\overline{B(0,e^A)} \times \overline{B(0,e^B)} \times
\overline{B(0,e^C)} = \{\,(x,y,z) : |x|\leq e^A, |y| \leq e^B,
|z|\leq e^C\,\}$.
\end{cl}

\begin{proof}[Proof of Claim \ref{cl3}]
Suppose for contradiction that there is a complex zero $(x,y,z)$
of $g$ with $(|x|,|y|,|z|) \in [0,e^A] \times [0,e^B] \times
[0,e^C]$, but no real zeros in the same region. As the zero locus
of $g$ is closed, we may assume that there is no complex zero
$(x',y',z')$ with $|x'|<|x|$, $|y'|<|y|$ and $|z'|<|z|$.  So the
power series of \[G(x,y,z) = \frac{1}{\frac{3}{2}(1-z)g(x,y,z)}\]
converges on $B(0,|x|) \times B(0,|y|) \times B(0,|z|)$ (we used
$C<0$ to conclude that the $(1-z)$ term doesn't vanish) and blows
up to $\infty$ as we approach $(x,y,z)$.

But the coefficients of $G(x,y,z)$ are all positive as they are
probabilities. So the series must also blow up as we approach
$(|x|, |y|, |z|)$ and thus $g(|x|,|y|,|z|)=0$ contradicting our
assumption that there are no zeroes in $[0,e^A] \times [0,e^B]
\times [0,e^C]$.
\end{proof}

\begin{proof}[Proof of Lemma]
We now apply the claims to the edge probability generating
function:
\begin{eqnarray}
G(x,y,z) & = & \sum_{i,j,k\geq
0}p(-i,-j,-k)x^{i}y^{j}z^{k}\nonumber\\
& = &
\frac{z^2}{\frac{3}{2}(1-z)(1+xyz-\frac{1}{3}(x+y+z+xy+xz+yz))}\nonumber\\
& = & \frac{z^2}{\frac{3}{2}(1-z)g(x,y,z)}.\nonumber
\end{eqnarray}
By Claim \ref{cl1} we only need to show that we can choose real
numbers $A,B,C$ so that
\begin{itemize}
\item $Ax_{0}+By_{0}+Cz_{0} < 0$
\item Both $1-z$ and $g(x,y,z)$ are not equal to zero for any $(x,y,z) \in
\{\,(x,y,z)\in \mathbb{C}^{3}:  |x|\leq e^{A}, |y|\leq e^{B},
|z|\leq e^{C} \,\}$.
\end{itemize}

We will now show that, for $(A_0,B_0,C_0)$ as in Claim \ref{cl2}
and $t$ positive and sufficiently small, $(A,B,C)=t(A_0,B_0,C_0)$
has the desired properties. We have $A x_0+B y_0+C z_0<0$ and
$C<0$ because the analogous properties hold for $(A_0,B_0,C_0)$.
All that remains is to show that $g(x,y,z)$ does not vanish for
$(x,y,z) \in \overline{B(0,e^A)} \times \overline{B(0,e^B)} \times
\overline{B(0,e^C)}$. By Claim \ref{cl3} it is enough to show that
$g$ has no zeroes on $[0,e^A] \times [0,e^B] \times [0,e^C]$. The
identity
$$g(x,y,z)=\frac{(1-x)(1-yz)+(1-y)(1-xz)+(1-z)(1-xy)}{3}$$
shows that $g(x,y,z) \neq 0$ on $[0,1)^3$. Writing $x=e^\alpha$,
$y=e^\beta$, $z=e^{\gamma}$, we have $g(x,y,z)=\alpha \beta+\alpha
\gamma + \beta \gamma + O(\alpha+\beta+\gamma)^3$. So, near
$(1,1,1)$, the zero locus of $g$ looks like the cone $\alpha
\beta+\alpha \gamma + \beta \gamma=0$.

Let $\mathcal{L} = \{\,(x,y,z) : g(x,y,z)=0 \,\}$ be the zero
locus. We want to show that there is a number $t > 0$ such that
the point $(e^{A},e^{B},e^{C}) = (e^{tA_0},e^{tB_0},e^{tC_0})$ is
inside of $\mathcal{L}$. We can write
\begin{eqnarray*}
g(e^{tA_0},e^{tB_0},e^{tC_0}) & = & t(A_0 B_0+A_0 C_0 + B_0 C_0) +
t^{3}O(A_0 + B_0 + C_0)^3\\
 & = & t(A_0 B_0+A_0 C_0 + B_0 C_0) + t^{3}O(-1)\\
 & = & t(A_0 B_0+A_0 C_0 + B_0 C_0 + t^{2}O(-1)).
\end{eqnarray*}
For fixed $(A_0,B_0,C_0)$, we can certainly choose a $t>0$ small
enough to guarantee that $A_0 B_0+A_0 C_0 + B_0 C_0 + t^{2}O(-1) >
0$. Therefore we have $(e^A,e^B,e^C) \notin \mathcal{L}$ and
moreover it is on the `inside' of $\mathcal{L}$ in that it is
closer to the origin than the locus.

Now we'd like to say $g$ will not vanish on $[0,e^A] \times
[0,e^B] \times [0,e^C]$. Since it is nonzero on $[0,1)^3$, we just
need to check that $g$ is not zero on $[1,e^A]\times [0,e^B]\times
[0,e^C]$, where we can take $B,C < 0$. Let $x = x(y,z)$ be a
parameterization of the zero locus. Then we have
\[x(y,z) = \frac{1/3(y+z+yz) - 1}{yz-1/3(1+y+z)}.\] It will be
enough to show that for any $(y,z)\in [0,e^B]\times[0,e^C]$, \[
e^A < x(e^B,e^C) \leq x(y,z).\] Let $(y,z)$ be any pair such that
$0\leq y\leq e^B$, $0\leq z\leq e^C$ and $e^A < x(y,z)$. The pair
$(e^B,e^C)$ is such an example. Then we will show that
$x(y,z-\epsilon)
> x(y,z)$ for any $0< \epsilon \leq z$, and similarly,
$x(y-\epsilon, z) > x(y,z)$ for any $0< \epsilon \leq y$. We have
\begin{eqnarray*}
x(y,z-\epsilon) & = & \frac{1/3(y+z-\epsilon+y(z-\epsilon)) -
1}{y(z-\epsilon)-1/3(1+y+z-\epsilon)}\\
 & = & \frac{1/3(y+z+yz) - 1 - 1/3\epsilon(y+1)}{yz-1/3(1+y+z) -
 \epsilon(y-1/3)}.
\end{eqnarray*}
Because both $y$ and $z$ are less than one, the numerator and
denominator of $x(y,z)$ are both negative. It now suffices to
check that $\epsilon(y-1/3) < 1/3\epsilon(y+1)$. This amounts to
saying $y < 1$, which is true by supposition. Hence,
$x(y,z-\epsilon) > x(y,z)$. By symmetry, $x(y-\epsilon, z)>
x(y,z)$ as well, and the lemma is proved.
\end{proof}

The asymptotics of multivariate generating functions is described in
some generality in the sequence of papers \cite{PW1}, \cite{PW2},
\cite{PW3}, by Robin Pemantle and Mark Wilson. Unfortunately, their methods do not apply to our generating function. Refer to the complex 2-manifold given by $1+xyz-1/3(x+y+z+xy+xz+yz) = 0$ as the \emph{singular variety}: what we refer to above as the zero locus. Our analysis requires studying the behavior of this surface near the point $(1,1,1)$. This point is a singularity of the variety that, as we have seen, resembles a cone point locally. The papers of Pemantle and Wilson give asymptotics near smooth points of the singular variety and near multiple points of the singular variety, but not near singular points of the singular variety as in our case. We hope that an extension of their
techniques can be used to prove a stronger version of Theorem \ref{thm:wac}.\footnote{The introduction to \cite{PW1} refers to the case of a singular point of the singular variety as a tractable problem to be handled in future work. Personal communication from Pemantle reveals that work on the case where the singular point locally resembles a quadratic surface, e.g. as in our case of the cone, may be finished very soon. See section \ref{sec:aztec} and mention of current work of Cohn and Pemantle.} In particular, we hope for a theorem that describes
the asymptotic probabilities in any location (rather than just outside of the circle), similar to Theorem 1 of
\cite{CEP}.

\section{Domino tilings of Aztec diamonds}\label{sec:aztec}

We now draw parallels between the examination of the behavior of
large groves on standard initial conditions, and the behavior of
tilings of large Aztec diamonds.  This approach yields no new
results for Aztec diamonds, but presents an alternative approach
to their study. In this section we derive a generating function for the probabilities $p_{n}(i,j)$ that position $(i,j)$ in a tiling of an Aztec diamond of order $n$ is covered by a particular
type of horizontal domino. Weak asymptotics for the function we
will derive are discussed as an example in \cite{PW1}, and ongoing work of Henry Cohn and Robin Pemantle seeks to give a full asymptotic expansion.\footnote{Thorough analysis of the generating function presented in this section requires analysis of the singular variety near a point that is locally a cone.} The first derivation of the function is due to James Propp and Alexandru Ionescu, though their (different) derivation has never been
published. Some recursive formulas for $p_{n}(i,j)$ are given in
\cite{propp}, and are the inspiration for our derivation of the
edge probabilities for groves.

\subsection{Tilings of Aztec diamonds and the octahedron
recurrence}

\begin{figure} [h]
\centering
\includegraphics[scale=.5]{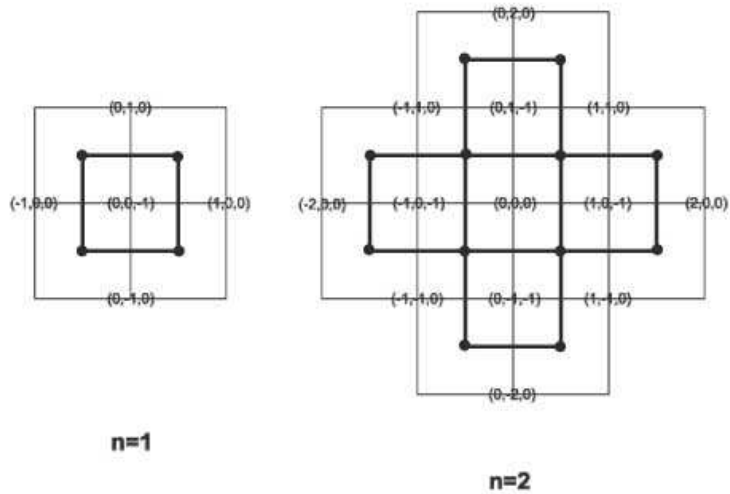}
\caption{Aztec diamond graphs of order $n$=1 and $n=2$.}
\label{fig:adgraph}
\end{figure}

We begin by describing precisely how tilings of Aztec diamonds are
encoded in the terms of polynomials generated by the octahedron
recurrence.\footnote{See Speyer's paper, \cite{Sp}, for more on encoding graphs with this recurrence.} First of all, rather than considering tilings of an
Aztec diamond, we prefer to consider perfect matchings of the dual
graph of this region. We call such a graph an \emph{Aztec diamond
graph}. The Aztec diamond graph of order $n$ has as its vertices
the set $\{\,(i\pm 1/2, j \pm 1/2)\, |\, i,j \in \mathbb{Z}, |i|+|j|
\leq n-1\,\}$.  In other words, it is the set of centers of the unit
squares that compose the Aztec diamond of order $n$.  Each vertex is
connected with an edge to its nearest horizontal and vertical
neighbors.  The \emph{faces} of the Aztec diamond graph are the
points $(i,j)$ such that $|i|+|j| \leq n$. For example, see Figure
\ref{fig:adgraph}.

\begin{figure} [h]
\centering
\includegraphics[scale=.4]{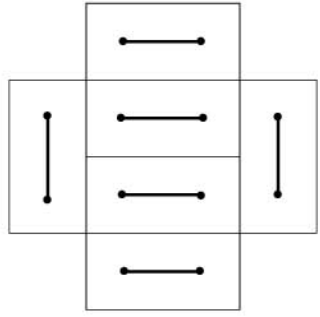}
\caption{Tiling of an Aztec diamond of order 2 and its dual
matching.} \label{fig:adex}
\end{figure}

Recall that the octahedron recurrence is:
\[g_{i,j,n}g_{i,j,n-2} = g_{i-1,j,n-1}g_{i+1,j,n-1} +
g_{i,j-1,n-1}g_{i,j+1,n-1}.\] We initialize the
octahedron recurrence by setting $g_{i,j,n} = x_{i,j,n}$ for $n =
0,-1$.  Then \[g_{0,0,n} = \displaystyle \sum_{\substack{\mbox{\tiny tilings
} T \\ \mbox{ \tiny of order } n }}\!\!m(T),\] where $m(T)$ is a Laurent
monomial in the variables $x_{i,j,\delta}$ where $\delta = 0$ if $i+j+n$
even, $\delta = -1$ if $i+j+n$ odd. Each monomial is of the form \[m(T) = \prod_{|i|+|j|\leq n}
x_{i,j,\delta}^{1-\deg(i,j)}\] where $\deg(i,j)\in \{0,1,2\}$ is
the number of edges surrounding face $(i,j)$ in the dual matching
to the tiling $T$. For example, in Figure \ref{fig:adex} we see
the tiling associated with the monomial
\[\frac{x_{-1,-1,0}x_{1,-1,0}x_{-1,1,0}x_{1,1,0}}{x_{0,-1,-1}x_{0,1,-1}x_{0,0,0}}.\]

The assertion that Aztec diamond tilings can be so encoded is just
a special case of Theorem 5.5 of Eric Kuo's paper on graphical
condensation \cite{K}. Here we take the weight of a tiling $T$ of
order $n$ to be $m(T)$ as above.

\subsection{Domino shuffling}

Domino shuffling is the name given to the local move that can be
used to generate random tilings of Aztec diamonds.  In a tiling of
an Aztec diamond there are four types of dominoes, which may be
characterized as follows.  Make a checkerboard coloring of the
Aztec diamond of order $n$ by making the leftmost square in each
row of the top half of the diamond be white.  A horizontal domino
is \emph{north-going} if its leftmost square is white,
\emph{south-going} if its leftmost square is black.  Similarly, a
vertical domino is \emph{east-going} if its topmost square is
black, \emph{west-going} if its topmost square is white.  Domino
shuffling is described in detail in \cite{propp}, but basically
each domino slides in the direction indicated by its name.
North-going dominoes take one step north, south-going dominoes
head south, and so on.  However, two dominoes may not pass through
one another.  If a north- and south- or east- and west-going pair
of dominoes collide, they annihilate one another, and the
resulting hole and any other holes opened by sliding dominoes is
filled in with a pair of dominoes: two horizontals with
probability 1/2, or two verticals with probability 1/2. Using this
description we have,
\begin{thm}\label{thm:ds}
The horizontal edge probabilities are given recursively by
$p_{n}(i,j) = p_{n-1}(i,j) + \frac{1}{2}E_{n-1}(i,j)$. Thus,
$\displaystyle p_{n}(i,j) =
\frac{1}{2}\sum_{l=1}^{n-1}E_{l}(i,j)$.
\end{thm}
The proof of the theorem can be found in \cite{propp}. It follows
more or less directly from the definition of domino shuffling,
where $E_{n}(i,j)$ is the \emph{net creation rate} (see
\cite{EKLP}, \cite{propp}).

\subsection{Another generating function}
By differentiating the uniformly weighted version of the
octahedron recurrence
\[g_{i,j,n+1}g_{i,j,n-1} =
\frac{1}{2}\left(g_{i-1,j,n}g_{i+1,j,n}+g_{i,j-1,n}g_{i,j+1,n}\right),\]
and because
\[E_{n}(i_{0},j_{0}) = \frac{\partial}{\partial x}\Big(g_{0,0,n}\Big)\Big|_{x_{i,j} =
1}\] we obtain
\begin{eqnarray}
E_{n+1}(i,j) + E_{n-1}(i,j) & = & \frac{1}{2}\left(E_{n}(i-1,j) +
E_{n}(i+1,j)\right) + \nonumber\\
& & \frac{1}{2}\left(E_{n}(i,j-1) + E_{n}(i,j+1)\right).\nonumber
\end{eqnarray}
From this recurrence and Theorem \ref{thm:ds} we get the
generating function:
\begin{eqnarray}
G(x,y,z) & = & \sum_{n\geq 0}\sum_{|i|+|j|\leq n}
p_{n}(i,j)x^{i}y^{j}z^{n}\nonumber \\
 & = &
 \frac{z/2}{(1-yz)(1+z^2-\frac{z}{2}(x+x^{-1}+y+y^{-1}))}.\nonumber
\end{eqnarray}

This is the form of the generating function used as an example in
\cite{PW1}. A weak arctic circle theorem like our Theorem
\ref{thm:wac} follows directly from that example. Probabilities
throughout the diamond could be extracted from this function in
principle, and current work of Cohn and Pemantle seeks to carry out this more difficult analysis.

\section{Biased groves, or, groves with a drift}
\begin{figure} [h]
\centering
\includegraphics[scale=.4]{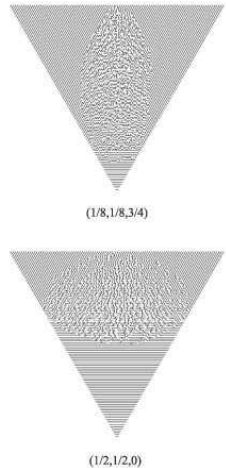}
\caption{Groves with bias $(\alpha,\beta,\gamma)$.}
\label{fig:bias}
\end{figure}

Another variation on the model presented in this paper would be to
study ``groves with a drift," i.e., rather than having the random
choice in grove shuffling be uniform, we make a biased
choice.\footnote{See \cite{CEP} for analysis of the same variation
in Aztec diamond tilings.} For example, say we choose the two diagonal edges
with probability $\alpha$, the horizontal and one of the diagonal edges with probability $\beta$,
and the horizontal and the other diagonal edge with probability $\gamma = 1-\alpha - \beta$. This
bias is reflected in the cube recurrence by setting
\[f_{i,j,k}f_{i-1,j-1,k-1} = \alpha f_{i-1,j,k}f_{i,j-1,k-1} +
\beta f_{i,j-1,k}f_{i-1,j,k-1} + \gamma
f_{i,j,k-1}f_{i-1,j-1,k}.\] The generating function for biased
creation rates is then
\begin{eqnarray}
F_{\alpha,\beta,\gamma}(x,y,z) & = & \sum_{i,j,k \geq 0}
E_{\alpha,\beta,\gamma}(-i,-j,-k)x^{i}y^{j}z^{k}
\nonumber \\
& = &
\frac{1}{1+xyz-\alpha(x+yz)-\beta(y+xz)-\gamma(z+xy)}.\nonumber
\end{eqnarray}
 We can also derive the generating function for biased horizontal edge probabilities:
\begin{eqnarray*}
G_{\alpha,\beta,\gamma}(x,y,z) & = & \sum_{i,j,k \geq 0}
p_{\alpha,\beta,\gamma}(-i,-j,-k)x^{i}y^{j}z^{k} \\
& = &
\frac{(\beta+\gamma)z^2 F_{\alpha,\beta,\gamma}(x,y,z)}{1-z}.
\end{eqnarray*}
As seen in Figure \ref{fig:bias}, the arctic circle is just a special case of what one might call the ``arctic ellipse,'' dependent on the free parameters $\alpha, \beta$. Using the same approach as in the unbiased case where $\alpha = \beta = 1/3$, we can prove the following theorem without too much difficulty.
\begin{thm}\label{thm:wac2}
The boundary of the frozen region for (rescaled) groves with a drift is given
by the intersection of the plane $x+y+z = -1$ with the surface
\[rs + rt + st = \frac{r^2 + s^2 + t^2}{2},\] where
$r = (\beta+\gamma)x$, $s =(\alpha+\gamma)y$, and $t=(\alpha+\beta)z$.
\end{thm}

Aside from simply describing the shape of the frozen region in the biased situation, it may also be useful to observe how the area of the temperate zone varies in $\alpha$ and $\beta$ (with respect to the area of the entire grove). Specifically, we can compute the ratio of the area of the temperate zone, $A_{\circ}$, to the area of the entire grove, $A_{\nabla}$, as a function of $\alpha$ and $\beta$. Given that $\gamma = 1-\alpha-\beta$,  \[\rho(\alpha,\beta) = \frac{A_{\circ}}{A_{\nabla}}= \frac{\pi(\alpha + \beta)(\alpha + \gamma)(\beta + \gamma)}{\left( (\alpha+\beta)(\alpha + \gamma) + (\alpha + \beta)(\beta + \gamma) + (\alpha + \gamma)(\beta + \gamma) \right)^{3/2}}.\]

\section{Speculation on statistics of groves}

As mentioned, we hope to apply the methods of Pemantle and Wilson
to determine asymptotic probabilities throughout a random grove.
Based on computer experiments and the similarity of groves and
Aztec diamond tilings seen so far, we believe a formula for such
probabilities exists. Another future aim is to apply the methods
of growth models and statistical mechanics to groves, in the style
of Johansson \cite{J1}, \cite{J2}, or more recently Kenyon,
Okounkov, and Sheffield \cite{Ke}, \cite{KO}, \cite{KOS},
\cite{KS}.

\subsection{Randomly growing Young diagrams}
Perhaps something can be proved about the variance of the boundary
circle by interpreting groves in a more familiar setting. One
clever way for determining the boundary of the frozen region for
Aztec diamond tilings is to look at a frozen corner as a randomly
growing Young diagram. See \cite{JPS}, \cite{J1}, \cite{J2} for
this interpretation. A nearly identical projection of the frozen
region of a grove yields some sort of randomly growing Young
diagram, but it seems to follow more intricate rules of growth
than those of Aztec diamond tilings.

\begin{figure} [h]
\centering
\includegraphics[scale=.7]{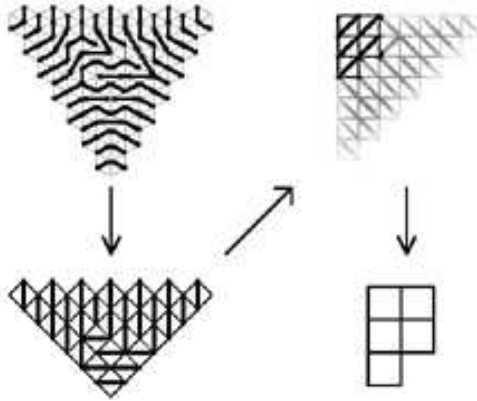}
\caption{Projecting a frozen corner to a Young diagram.}
\label{fig:proj}
\end{figure}

If we project the grove onto the plane $\mathbb{R}^{3}/(0,0,1)$ we
see a triangular array of boxes with two types of diagonal edges.
Let us put the corner box at the upper-left and index it as (0,0).
Then we index the rows by $0\leq i \leq n-1$, the columns by
$0\leq j \leq n-1$, so that each box has index $(i,j)$ with $i+j
\leq n-1$. A diagonal edge from the top left corner to the bottom
right corner of a box corresponds to a short edge in the grove.
Edges from the top right to bottom left correspond to long edges
in the grove. The box $(i,j)$ is \emph{frozen} if it contains a
long edge and all the boxes $(i',j')$ contain long edges, $0\leq
i' \leq i$, $0\leq j'\leq j$. Clearly the collection of all frozen
boxes is a Young diagram. We would like to be able to describe how
this Young diagram grows under grove shuffling.

\begin{figure} [h]
\centering
\includegraphics[scale=.8]{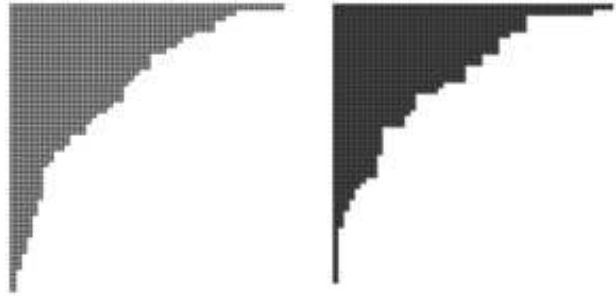}
\caption{On the left, the Young diagram corresponding to a random
grove of order 100. On the right, the Young diagram corresponding
to a random tiling of an Aztec diamond of order 100.}
\label{fig:gyd}
\end{figure}

In a randomly growing Young diagram, we call box $(i,j)$ a
\emph{growth position} if boxes $(i-1,j)$ and $(i,j-1)$ are both
frozen and $(i,j)$ is not frozen (we use the convention that boxes
of the form $(i,-1), (-1,j)$ are always frozen). In the case of
Aztec diamond tilings, the growth of a frozen corner under domino
shuffling corresponds exactly to adding a new box at each growth
position independently with probability 1/2. With groves this is
not the case, though it may not be clear from Figure
\ref{fig:gyd}. In fact, two groves may have the same Young diagram
projection, yet grow very differently upon shuffling. It can
happen that, at a particular growth position, one grove will not
permit the addition of a new box under any circumstances, whereas
the other will add a box with probability 2/3.

Another difference seems to suggest the need for a new definition of growth position.
It is possible that a grove can project to a Young diagram with
growth position $(i,j)$, but after just one iteration of grove shuffling the projection adds a box not only to position $(i,j)$, but also position $(i+1,j)$. In fact for any $k$ there is
a grove of some (perhaps large) order $n$ with growth position
$(i,j)$, so that with positive probability its projection to a
Young diagram adds boxes $(i,j), (i+1,j), \ldots, (i+k,j)$ upon
one iteration of grove shuffling. Perhaps such situations are
outliers, but there is still much work to do in this direction.

\subsection{The nexus} 
\begin{figure} [h]
\centering
\includegraphics[scale=.6]{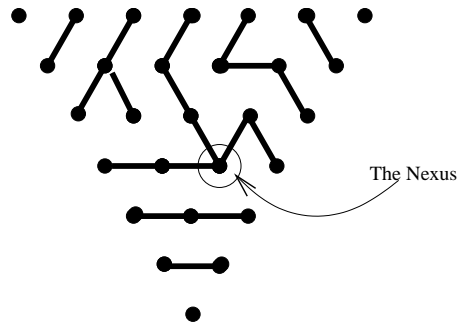}
\caption{The nexus of a grove.}
\label{fig:nex}
\end{figure}

Much of the statistical study of groves
is motivated by analogy with statistics for domino tilings of
Aztec diamonds. But there is at least one interesting feature of groves
that seems to have no analogy in the realm of Aztec diamonds. We conclude the paper with some observations about a unique vertex that is present in every grove. We call this special vertex the \emph{nexus}. Loosely speaking, the nexus is the vertex at the ``middle" of the unique tree connected to all three sides of the initial conditions. If the nexus and its incident edges are removed from the grove, then all three sides become disconnected from one another. In Figure \ref{fig:nex} we see the nexus highlighted.

\begin{defn}
Let $g$ be a grove of order $2n$ or $2n+1$, $n =1,2,\ldots$. The
\emph{nexus} of $g$ is the (even) vertex $v$ which is connected to
each of the three \emph{midpoints}: $(-n,-n,0)$, $(-n,0,-n)$, and
$(0,-n,-n)$, and for which each midpoint lies on a distinct branch
of the tree rooted at $v$.
\end{defn}

We would like to understand how the nexus moves during grove shuffling. The nexus takes some kind of random walk in the initial conditions, but it is not a simple random walk. In some sense the nexus is a ``stuttering'' random walker. In some situations the nexus takes a step in one of three directions with equal probability, in others it does not move at all, and in still others it moves deterministically. Specifically, if the nexus is at an ``up" vertex (see section \ref{sec:init}) then after one iteration of grove shuffling it will take a step in one of three directions, each with equal probability. If the nexus is at a ``down" vertex, then after one iteration of grove shuffling the nexus will not move (and so becomes a ``flat" vertex in the new grove). Strangely, if the nexus is at a flat vertex then depending on which edges lead from the nexus to the midpoints, the nexus will either move deterministically one step or it will not move. Ultimately we want to say something about the distance of the nexus from the center of the grove as the size of the grove gets very large.  How far from home does the nexus roam?

\section{Acknowledgements}
We would like to thank Jim Propp for bringing us together to work on this problem, and Robin Pemantle for his advice and encouragement.

\end{document}